\documentclass[letterpaper, 10 pt, conference]{ieeeconf}  

\IEEEoverridecommandlockouts                              
\usepackage{braket,amsfonts}

\usepackage{array}

\usepackage{pgfplots}

\usepackage{algorithmicx}
\usepackage[ruled]{algorithm}
\usepackage[noend]{algpseudocode}
\usepackage{dsfont}
\usepackage{graphicx,epstopdf}

\newtheorem{thrm}{Theorem}[section]
\newtheorem{lmm}[thrm]{Lemma}
\newtheorem{crllr}[thrm]{Corollary}
\newtheorem{proposition}[thrm]{Proposition}

\newtheorem{assumption}[thrm]{Assumption}

\newtheorem{xmpl}[thrm]{Example}

\usepackage{graphics} 
\usepackage{dsfont}
\usepackage{epsfig} 
\usepackage{mathptmx} 
\usepackage{times} 
\usepackage{amsmath} 
\usepackage{amssymb}  
\usepackage{mathtools}
\usepackage{siunitx}
\usepackage{subcaption}
\usepackage{tikz}
\usetikzlibrary{arrows, automata, backgrounds, calc, decorations.pathmorphing, decorations.shapes, fit, intersections, matrix, petri, plotmarks, positioning, shapes, snakes, through, trees}
\usepackage{xspace}
\usepackage{bold-extra}
\usepackage{bbm}
\usepackage[most]{tcolorbox}

\usepackage[skip=1pt]{caption}

\colorlet{texcscolor}{blue!50!black}
\colorlet{texemcolor}{red!70!black}
\colorlet{texpreamble}{red!70!black}
\colorlet{codebackground}{black!25!white!25}

\renewcommand{\baselinestretch}{0.96}

\tikzset{%
    state/.style={circle,fill=none,draw=black,text=black, minimum size = 29pt},
    const/.style={circle,fill=red!40,draw=black,text=black, minimum size = 29pt},
    triang/.style={fill=none,draw,inner sep=0pt,label distance=2pt,minimum size=2pt,circle}
}

\title{\LARGE \bf Identification of Cyclists' Route Choice Criteria
}

\author{Stefano Ardizzoni, Mattia Laurini, Rafael Praxedes, Luca Consolini, and Marco Locatelli  
\thanks{All authors are with the University of Parma, Department of Engineering and Architecture, Parco Area delle Scienze 181/A, 43124 Parma, Italy. E-mails:
{\tt\footnotesize \{stefano.ardizzoni, mattia.laurini, rafael.praxedes, luca.consolini, marco.locatelli\}@unipr.it}.
}%
\thanks{Project funded under the National Recovery and Resilience Plan (NRRP), Mission 4 Component 2
Investment 1.5 -- Call for tender No. 3277 of 30/12/2021 of Italian Ministry of University and Research
funded by the European Union -- NextGenerationEU.}%
\thanks{Award Number: Project code ECS00000033, Concession Decree No. 1052 of 23/06/2022 adopted by the
Italian Ministry of University and Research, CUP D93C22000460001,
``Ecosystem for Sustainable Transition in Emilia-Romagna'' (Ecosister).}}%
\date{}

\begin{document}

\maketitle

\begin{abstract}
The behavior of cyclists when choosing the path to follow along a road network is not uniform.
Some of them are mostly interested in minimizing the travelled distance, but some others may also take into account other features such as safety of the roads or pollution.
Individuating the different groups of users, estimating the numerical consistency of each of these groups, and reporting the weights assigned by each group to different characteristics of the road network, is quite relevant.
Indeed, when decision makers need to assign some budget for infrastructural interventions, they need to know the impact of their decisions, and this is strictly related to the way users perceive different features of the road network.
In this paper, we propose an optimization approach to detect the weights assigned to different road features by various user groups, leveraging knowledge of the true paths followed by them, accessible, for example, through data collected by bike-sharing services.
\end{abstract}

\section{Introduction}
The transition to more sustainable and green forms of transportation is increasingly becoming a priority in developed countries.
In particular, bicycles, electric bicycles, and electric scooters are a convenient mode of transportation for short range and urban travels.
Understanding how users (or simply cyclists, in what follows) choose their routes depending on various road characteristics is fundamental if one aims at increasing the number of cyclists (and, hence, decreasing the number of motor vehicles users) or improving the existing cycling infrastructure, helping decision makers take more informed actions when assigning budget for infrastructural interventions (see, e.g.,~\cite{LIU2019, ZUO2019}).

\subsection{Literature review}
A quantitative method for assessing the quality of roads from a cyclist's point of view is given by the concept of Bicycle Level of Service (BLOS), which has first been introduced in the late '80s -- early '90s in~\cite{DAVIS1987,EPPERSON1994}.
Its aim is to measure quantitatively several qualitative aspects of road segments with respect to cyclists' perception.
As we know from several studies (see, for instance,~\cite{BUEHLER2016,PUCHER2008}), cyclists may not use distance as the only objective function when choosing which route to follow.
For instance, the presence of bike facilities may heavily influence the route choice (see, e.g.,~\cite{HOWARD2001}), resulting in longer paths in which the amount of road sections with bike facilities seems to be maximized.
However, this is just one example of how the features of a road portion may influence the cyclists' choice.
In its original formulation, despite its innovative aspect, the BLOS was affected by some shortcomings like the lack of statistical calibration and a subjective methodology in assigning road features values.
Over the years, similar concepts have been developed such as that of Bicycle Compatibility Index (BCI) (see~\cite{HARKEY1998}), which aims at evaluating the suitability of roadways for accommodating both motor vehicles and bicycles, or that of bikeability (see, e.g.,~\cite{NACTO2014,WINTERS2013}), which has the goal of assessing how promotive an environment is for biking.
The same concept of BLOS has been further explored and studied by the scientific community, including in the BLOS formulation more and more aspects that may affect the perception and choice of the users of a bike network.
In early formulations, only the infrastructure aspects of road segments were considered together with bicycle flow interruptions.
Now, research is focusing on including exogenous factors in the BLOS or features on which decision makers cannot apply direct interventions, such as climate factors, presence of pollens, topographic features, but also pollution, noise, and so on (see, for instance,~\cite{KAZEMZADEH2020}).
One of the most critical aspects in BLOS is that of determining the weighting factors multiplying the quantities associated to the various considered aspects of a road section.
Obtaining a ``good'' set of coefficients can require data collection, surveying users, normalization and homogenization of different measurement scales, and it also requires validation and continuous calibration of the obtained formula.

\subsection{Statement of contribution}
In our work we assume that $r$ basic objective functions (i.e., road features) are given and that users consider a combination of such functions in order to determine the path to be followed.
This is equivalent to defining a BLOS formula in which only $r$ factors are involved.
The road network is represented by a graph and each basic objective function is defined assigning costs to all arcs of this graph.
We consider a BLOS formula which is a convex combination of the considered features, hence the coefficients of the BLOS formulas, also called weights in what follows, are all assumed to be in $[0, 1]$ and their sum is equal to one.
We assume that each user has their own $r$-dimensional weight vector, and follows a shortest path (SP) over the graph representing the road network, where the costs of the arcs are a convex combination of the $r$ basic costs, with coefficients of the convex combination corresponding to the entries of the weight vector. We also assume that users may have different behaviors and, thus, select their paths according to different weight vectors. Therefore, 
the goal of this paper is that of identifying both the set of weighting factors that users perceive, and the probability with which users would consider such weights.
The identification is based solely on traffic flow observations on (a subset of) arcs of the network.
In other words, we assume that there is not a unique BLOS formula that suits all users but we aim at identifying different BLOS formulas for different users segments.

Note that the graph considered in this work shares some similarities with the one presented in~\cite{STEINACKER2022}, in which a preference graph is considered.
There, the weight of each edge depends on a combination of various factors, which however are assumed to be known or measurable. In our work, we aim at estimating such quantities.
Other works use GPS data in order to determine the route-choice models such as~\cite{HOOD2011, WARGELIN2012}.
Others aim at optimizing the BLOS along paths used by cyclists of the network (see, for instance,~\cite{SMITH2012, DUTHIE2014}).
The present work serves as a preliminary method for BLOS identification which can be particularly useful when addressing the problem of bike network optimization, in which one wishes to maximize the benefit of infrastructural interventions given a limited budget.




\subsection{Paper organization}
The paper is structured as follows. In Section~\ref{sec:formul}, we formalize the problem.
More precisely, in Section~\ref{sec:known}, we consider a simplified version where the set of possible weights is assumed to be known in advance.
For this problem, we propose a bilevel optimization formulation, and derive a polynomial-time algorithm for its solution.
Next, in Section~\ref{sec:unknown}, we present the optimization problem with unknown set of weights, and we discuss some properties of the function to be minimized.
In Section~\ref{sec:data}, we discuss how the data needed for the problem definition can be collected, also pointing out possible difficulties and limitations of the proposed approach.
In Section~\ref{sec:algo}, we propose a solution algorithm for the problem presented in Section~\ref{sec:unknown}.
Finally, in Section~\ref{sec:compexp}, we present some preliminary experiments on synthetic data.

\section{Problem formulation}
\label{sec:formul}
We represent a bicycle network with a directed graph $G=(V,A)$. We denote the number of nodes and directed arcs by $n=|V|$ and $m=|A|$, respectively. The arc set  represents the roads used by cyclists, and the node set the intersections.
The network is assumed to be strongly connected.

Together with the network, we are also given a set $W\subset V \times V$, made up of origin-destination (O-D) pairs within the network. We associate to each pair $w=(o_w,d_w)\in W$ a demand value $u_w\in \mathbb{N}$, corresponding to the number of users that move from node $o_w$ to node $d_w$ traveling within the network at a given time of the day. Values $u_w$, with $w\in W$, may be known in advance but, in some cases, there might be the need to estimate at least some of them.
We denote by $\Pi_w$ the set of elementary directed paths from origin node $o_w$ to destination node $d_w$. To each $\pi \in \Pi_w$, we associate subset $A_{\pi}^w \subseteq A$ of arcs belonging to directed path $\pi$.           We associate to each arc $(i,j)\in A$ a flow $x_{ij}^w$ corresponding to the amount of users associated to pair $w$ traveling along the arc.
For all $(i,j) \in A$, the total flow along arc $(i,j)$ is denoted by 
\begin{equation}
	\label{eq:const1}
	x_{ij}=\sum_{w\in W} x_{ij}^w.
\end{equation}
Thus, we have the following vectors:
\begin{itemize}
	\item $x\in \mathbb{R}_+^{|A|}$, the vector whose components are the total flows $x_{ij}$ along the arcs of the network;
	\item $x^w\in \mathbb{R}_+^{|A|}$, $w\in W$, the vectors whose components are flows $x_{ij}^w$ of users associated to O-D pair $w\in W$ along the arcs of the network.
\end{itemize}
Moreover, we associate to each arc $(i,j)\in A$ a set of costs $c_{ij}= (c_1^{ij},\ldots,c_r^{ij})$, that represents the characteristics  of that road, such as the length, the security, environmental conditions, and so on. These $r$ costs will be called in what follows {\em basic} costs. The total cost of an arc is a convex combination of these values. We assume that this combination depends on the single user. This is because each cyclist can choose the best route differently, giving more importance to one feature rather than another. Note that, in order to take into account that different features have different units of measure and different magnitudes, we normalize the basic costs in such a way that
$\sum_{(i,j)\in A} c^{ij}_{h_1}=\sum_{(i,j)\in A} c^{ij}_{h_2}$ for all $h_1,h_2\in \{1,\ldots,r\}$.

We first consider a simplified problem where the set of convex combinations returning the arc costs for each user are known in advance. Later on, we will address the problem where such set is also to be computed.
\subsection{The case with known set of convex combinations/weights}
\label{sec:known}
The coefficients of a convex combination will be called in what follows {\em weights}. In this subsection
we assume that the set of feasible weights $P= \{p^\ell = (p_1^\ell,\ldots,p_r^\ell)\geq 0 \mid (\forall \ell \in \{1,\ldots,q\})\ \sum_{h=1}^r p_h^\ell=1\}$ is known in advance. Therefore, the total cost of arc $(i,j)$ for a user that chooses weight $p^\ell$ is 
\[c_{ij}^\top p^\ell = \sum_{h=1}^{r} c_h^{ij} p_h^\ell. \]

We assume that the traffic is not congested. Hence, all vehicles moving from $o_w$ to $d_w$ follow the SP in $G$. However, the SP depends on the chosen weights, so it may be different for each user.

We assume that users choose between feasible weights according to a certain probability distribution. Therefore, to each weight $p^\ell$ we associate a value $\alpha_\ell$, which represents the probability that any user chooses that particular convex combination to calculate the SP. 

Therefore, $P$  represents the set of possible cyclists' route choice criteria and $ \{\alpha_\ell\}_{\ell \in \{1,\ldots,q\}}$ is the probability distribution that describes how many users choose them.

To each $p^\ell \in P$, with $\ell \in \{1,\ldots, q\}$, and each $w\in W$, we associate a vector $x^{w,\ell}\in \mathbb{R}^{|A|}$, whose components are flows $x_{ij}^{w,\ell}$ along the arcs of the network of users associated to O-D pair $w \in W$ whose selected weight vector is $p^\ell$.
For all $(i,j)\in A$, and $w \in W$, we have the following constraints, linking variables $x_{ij}^{w,\ell}$ and $x_{ij}^w$,

\begin{equation}
	\label{eq:const2}
	x_{ij}^w=\sum\limits_{\ell =1}^{q}  x_{ij}^{w,\ell}.
\end{equation}

For each pair $w\in W$, for each node $i \in V$, and for each weights combination $p^\ell \in P$, the following flow conservation constraints hold:
\begin{equation}
	\label{eq:const3}
	\sum_{(i,j)\in A} x_{ij}^{w,\ell} - \sum_{(j,i)\in A} x_{ji}^{w,\ell}= \alpha_\ell \, q_i^w,
\end{equation}
where
\[
q_i^w :=
\begin{cases}
u_w	& i=o_w \\
-u_w	& i=d_w \\
0	& i\in V\setminus \{o_w,d_w\}. 
\end{cases}
\]

We further impose the non-negativity constraints
\begin{equation}
\label{eq:const4}
(\forall w \in W)\ (\forall \ell \in \{1,\ldots, q\})\ x, x^w,x^{w,\ell}\geq 0.
\end{equation}

Constraints~\eqref{eq:const3} can be written in the matrix form
\[
N x^{w,\ell}=\alpha_\ell q^w,
\]
where $N$ is the node-arc incidence matrix of the graph, and $q^w$ is the vector of the demand related to pair $w \in W$.

We assume that we do not measure all flows, but only a subset of them $\bar{A} \subset A$. We denote with $\bar{x}_{ij}, (i,j)\in \bar{A}$ the measured flows. Our aim is to estimate $\{\alpha_\ell\}_{\ell \in \{1,\ldots,q\}}$  from the knowledge of  $\bar{x}_{ij}, (i,j)\in \bar{A}$.

To this end, we minimize the squared distance between measured and estimated flows, so that we end up with the following bilevel optimization problem:


\begin{align}
    g(P)=\min_{\alpha} \sum_{(i,j) \in \bar{A}} \left[x_{ij} - \bar{x}_{ij}\right]^2 \label{eq:FO_up}
\end{align}
\quad \quad \quad s.t.
\begin{align}
    x_{ij} = \sum_{\ell = 1}^{q} \sum_{w \in W} x_{ij}^{w,\ell} && (i,j) \in \bar{A} \label{eq:constr1_up} \\
    \sum_{\ell = 1}^{q} \alpha_\ell = 1 \label{eq:constr2_up} \\
    \alpha_\ell \geq 0 && \ell \in \{1,\ldots,q\} \label{eq:constr4_up} \\
    x^{w,\ell} \in S(\alpha_{1},\ldots,\alpha_{q}) && \ell \in \{1,\ldots,q\}, w \in W \label{eq:constr5_up}
\end{align}
\begin{align}
    S(\alpha_{1},\ldots,\alpha_{q}) = \arg\min \sum_{\ell = 1}^{q} \sum_{w \in W} \sum_{(i,j) \in A} \left( c_{ij}^\top p^\ell \right) x_{ij}^{w,\ell} \label{eq:FO_lw}
\end{align}
\quad \quad \quad s.t.
\begin{align}
    N x^{w,\ell} = \alpha_\ell q^w && \ell \in \{1,\ldots,q\}, w \in W \label{eq:constr1_lw} \\
    x^{w,\ell} \geq 0 && \ell \in \{1,\ldots,q\}, w \in W. \label{eq:constr2_lw}
\end{align}
Note that the optimal value $g(P)$ of the bilevel problem depends on the set of weights $P$, assumed to be known in advance.
For the upper-level model, objective function~\eqref{eq:FO_up} aims to minimize the error between the calculated and measured flows. Constraints~\eqref{eq:constr1_up} compute the total arc flows, while constraints~\eqref{eq:constr2_up} and~\eqref{eq:constr4_up} impose that $\alpha$ is a probability distribution over the set $P$ of weights. 
Constraints~\eqref{eq:constr5_up} impose that the estimated flows must be optimal with respect to the lower-level problem parameterized by the upper-level variables $\alpha_\ell$.

For the lower-level model, objective function~\eqref{eq:FO_lw} aims to minimize the total travel cost, subject to the fulfillment of the O-D pairs demands guaranteed by constraints~\eqref{eq:constr1_lw}. Constraints~\eqref{eq:constr2_lw} define the domain of the lower level variables.  

Problem~\eqref{eq:FO_up}--\eqref{eq:constr2_lw} can be solved quite efficiently. 
First we observe that the lower-level problem can be split into the following $q|W|$ subproblems: for each $\ell \in \{1,\ldots,q\}$ and each $w\in W$, solve:
\begin{align*}
\min	& \sum_{(i,j) \in A}  \left( c_{ij}^\top p^\ell \right)  x_{ij}^{w,\ell} \\
	& N x^{w,\ell} = \alpha_\ell q^w \\
	& x^{w,\ell} \geq 0.
\end{align*}
The solution of this problem is obtained by first detecting the SP from $o_w$ to $d_w$ with cost of each arc $(i,j)$ equal to $c_{ij}^\top p^\ell$. Once the SP has been detected, we send a flow equal to $\alpha_\ell u^w$ along the arcs of the path. More precisely, we proceed as follows.
Let $S^{w,\ell} \subset A$ be a SP from $o_w$ to $d_w$ based on the weight vector $p^\ell$.
Additionally, let
\begin{equation}
\label{eq:fij}
    f_{ij}^{w,\ell} = \begin{cases}
        u^w, \hspace{5 mm} (i,j) \in S^{w,\ell} \\
        0, \hspace{7 mm} \text{otherwise.}
    \end{cases}
\end{equation}
Next, we define a matrix $M_{|A| \times q}$ whose elements correspond to the sum of flows for each $w \in W$, considering the arc $(i,j)$ and the weight $p^\ell$, that is,
\begin{equation}
\label{eq:matrixM}
(\forall \ell \in \{1,\ldots,q\})\ (\forall (i,j) \in A)\ M_{(i,j),\ell} = \sum_{w \in W} f_{ij}^{w,\ell}.
\end{equation}
Then,
\begin{equation}
\label{eq:compxij}
(\forall (i,j)\in A)\ x_{ij}=\sum_{\ell=1}^q  M_{(i,j),\ell} \alpha_\ell,
\end{equation}
or, in matrix form $x=M\alpha$.
If we denote by $M_{\bar{A}}$ the submatrix obtained by considering only the rows of $M$ in $\bar{A}$, the upper-level problem reduces to:
\begin{equation}
\label{eq:convexQP}
\begin{aligned}
\min_{\alpha} & \|M_{\bar{A}}\alpha-\bar{x}\|^2 \\
& \sum_{\ell=1}^q \alpha_\ell=1 \\
& \alpha\geq 0,
\end{aligned}
\end{equation}
which is a convex quadratic problem, solvable by different available commercial solvers like, e.g., {\tt{Gurobi}}~\cite{gurobi}.
In summary, the algorithm to compute the optimal probability distribution over a fixed set $P$ of weights is the following:\\

\noindent
{$[\alpha_P,g(P)]\ =\ $\textsc{Identification}($P$)}
\begin{description}
\item[Step 1]  For each $w\in W$ and $p^\ell\in P$, compute the SP from $o_w$ to $d_w$ with cost
$c_{ij}^\top p^\ell$ associated to each arc $(i,j)$;
\item[Step 2 ]  Compute matrix $M\in \mathbb{Z}^{|A|\times |P|}$ with entries defined in~\eqref{eq:matrixM}, where values $f_{ij}^{w,\ell}$ are defined in~\eqref{eq:fij};
\item[Step 3 ] Solve convex Quadratic Programming (QP) Problem~\eqref{eq:convexQP}.
\end{description}
The overall complexity of this algorithm is stated in the following proposition.
\begin{proposition}
\label{prop:compl}
The complexity of the proposed algorithm for the case of known weights is 
\[
O\left(|W||P|\left(|A|+|V|\log|V|\right)+|P|^3 L\right),
\]
where $L$ is the bit size of the input of the convex QP.
\end{proposition}
\begin{proof}
Step 1 of the algorithm requires the solution of $|W||P|$ SP problems, so that the complexity of Step 1, if Dijkstra's algorithm is employed to solve the SP problems, is:
$O\left(|W||P|\left(|A|+|V|\log(|V|)\right)\right)$. Step 2 requires a time $O(|W||P||V|)$ since, for each $w\in W$ and each $p^\ell\in P$, only the entries of column $\ell$ of matrix $M$ associated to the arcs in the SP from $o_w$ to $d_w$ are updated, and the SP contains at most $|V|$ arcs.
Finally, the convex QP problem belongs to the class of problems for which, in~\cite{monteiro1989interior}, it is shown that the computing time for their solution is $O\left(|P|^3 L\right)$. 
\end{proof}
Note that this complexity result shows that for small $|P|$ values the major cost is represented by the solution of the SP problems, but as $|P|$ increases, the major cost becomes the solution of the convex QP problem.
\subsection{The case of unknown weights}
\label{sec:unknown}
If the set $P$ of weights is not known in advance, then a further optimization has to be performed, searching for a set $P$ with lowest possible value $g(P)$ (i.e., lowest possible distance between observed and estimated flows).
The value of $g$ can be reduced by: (i) enlarging the set of weights and/or (ii) perturbing the current weights. 
Enlarging the set of weights allows to reduce $g$ because of the monotonicity property of $g$, proved in the following proposition.\\

\begin{proposition}
\label{prop:monot}
Let $P'\supset P$. Then, $g(P')\leq g(P)$.
\end{proposition}
\begin{proof}
Let us denote by $M_{\bar{A},P}$ the restriction of matrix $M$ with rows in $\bar{A}$ and columns in $P$.
Moreover, let $\Delta_{|P|}=\left\{\alpha\in \mathbb{R}_+^{|P|} \mid \sum_{i=1}^{|P|} \alpha_i=1\right\}$. Then,
\begin{align*}
g(P) = & \min_{\alpha_P\in \Delta_{|P|}}  \|M_{\bar{A},P}\alpha_P-\bar{x}\|^2 \\
g(P') = & \min_{\alpha_{P'}\in \Delta_{|P'|}} \|M_{\bar{A},P'}\alpha_{P'}-\bar{x}\|^2.  
\end{align*}
Since $P\subset P'$, then we have $\alpha_{P'} = [\alpha_P, \alpha_{P' \setminus P}]$ and $M_{\bar{A},P'}= [M_{\bar{A},P}| M_{\bar{A},P' \setminus P'}]$.  Let $\bar{\alpha}_P$ be a feasible solution of the optimization problem~\eqref{eq:convexQP} with set of weights $P$. If we set $\bar{\alpha}_{P' \setminus P}=0$,  then $\bar{\alpha}_{P'} = [\bar{\alpha}_P, \bar{\alpha}_{P' \setminus P}]$ is a  feasible solution of the optimization problem~\eqref{eq:convexQP} with set of weights $P'$, with the same objective function value as $\bar{\alpha}_P$. Therefore, to each feasible solution of the first problem with set $P$, we can associate a feasible solution of the second problem with set $P'$, and the two solutions have the same objective function value. Then the inequality
$g(P')\leq g(P)$ immediately follows.
\end{proof}
In fact, a rather similar proof can be applied to reduce the number of weights.
\begin{crllr}
\label{coroll:1}
	Let $P$ be a set of weights and $\alpha$ be an optimal solution of the optimization problem~\eqref{eq:convexQP}. If $\bar{P}\subset P$ is such that, for all $p^\ell \in P\setminus \bar{P}$, $\alpha_\ell =0$, then $g(P)=g(\bar{P})$.
\end{crllr}
\begin{proof}
The optimal solution $\alpha$ can be written as $[\alpha_{\bar{P}}, \alpha_{P\setminus \bar{P}}]$, where, by assumption, $\alpha_{P\setminus \bar{P}}=0$. Then, $\alpha_{\bar{P}}$ is a feasible solution of~\eqref{eq:convexQP} with set of weights $\bar{P}$ and its objective function value is equal to $g(P)\leq g(\bar{P})$. Then, $\alpha_{\bar{P}}$ is also an optimal solution of~\eqref{eq:convexQP} with set of weights $\bar{P}$, and $g(P)= g(\bar{P})$ holds.
\end{proof}
According to Proposition~\ref{prop:monot}, we can reduce $g$ by expanding the set of weights. However, a large set of weights $P$ has at least two drawbacks. The first one is that the complexity result stated in Proposition~\ref{prop:compl} shows that the computing times for the algorithm calculating value $g(P)$ grow as $|P|^3$.
The second drawback is that, for the sake of interpretability, large $|P|$ values should be discouraged. Note that the result stated in Corollary~\ref{coroll:1} allows reducing the set of weights by discarding all weights with null probability.
\newline\newline\noindent
An alternative way to reduce $g$ is by keeping fixed the cardinality $q$ of $P$ and by perturbing the weights in $P$. Then, we can introduce a function 
\[
\bar{g}\ :\ \Delta_r^q\rightarrow \mathbb{R}_+,
\]
where:
\[
\Delta_r=\left\{p\in \mathbb{R}_+^r \mid \sum_{h=1}^r p_j=1\right\},
\]
that is, $\Delta_r$ is the $r$-dimensional unit simplex, defined as follows: if $P=\{p^1,\ldots,p^q\}$, where $p^\ell\in \Delta_r$, with $\ell \in \{1,\ldots,q\}$, then
$\bar{g}(p^1,\ldots,p^q)=g(P)$. 
Hence, the problem of identifying the best set of weights with fixed cardinality $q$ can be reformulates as follows:
\[
\min_{(p^1,\ldots,p^q)\in \Delta_r^q}\ \bar{g}(p^1,\ldots,p^q).
\]
Unfortunately, we cannot employ gradient-based methods even to detect local minimizers of $\bar{g}$. Indeed, 
we will show that this function is not continuous and is piecewise constant.
To this end, we first introduce an assumption.

\begin{assumption}
\label{ass:dist}
For each $w\in W$, recall that $\Pi_w$ is the finite collection of paths between $o_w$ and $d_w$. For some $\pi\in \Pi_w$, let $\lambda_h(\pi)=\sum_{(i,j)\in \pi} c^{ij}_h$, $h \in \{1,\ldots,r\}$, be the cost of $\pi$ with respect to the $h$-th basic cost. Then, we assume that for each  $w\in W$, there do not exist two distinct paths $\pi,\pi'\in \Pi_w$ such that
$\lambda_h(\pi)=\lambda_h(\pi')$ for all $h\in \{1,\ldots,r\}$.
\end{assumption}
The assumption simply states that there do not exist two distinct paths equivalent with respect to all the $r$ basic costs. If such two paths existed, they would be indistinguishable, since their cost would be the same for all possible weights $p\in \Delta_r$. 
Now, we can prove the following lemma.
\begin{lmm}
\label{lem:null}
For each $w\in W$ and $p^\ell\in \Delta_r$ let $S^{w,\ell}$ denote the set of SPs when the cost of each arc 
$(i,j)$ is $c_{ij}^\top p^\ell$. 
Under Assumption~\ref{ass:dist}, the set
\[
\Delta_r(w,\ell) = \left\{p^\ell\in \Delta_r \mid |S^{w,\ell}|>1\right\},
\]
has null measure in $\Delta_r$.
\end{lmm}
\begin{proof}
Let $p^\ell\in \Delta_r(w,\ell)$, then there exist $\pi_1$,$\pi_2$ $\in S^{w,\ell} $. It follows that these two paths have the same cost
\[
\sum_{h=1}^{r} p_h^\ell \lambda_h(\pi_1) = \sum_{h=1}^r p_h^\ell \lambda_h(\pi_2),
\]  
and so  $\sum_{h=1}^r p_h^\ell (\lambda_h(\pi_1) - \lambda_h(\pi_2)) =0$.  We want to prove that $\dim(\Delta_r(w,\ell)) < \dim(\Delta_r) = r-1$. Note that the space $\Delta_r(w,\ell)$ is the set of weights that satisfy the system of two linear equations 
\[
R p^\ell = \left[\begin{matrix}
1 \\
0
\end{matrix}\right],
\] 
where    
\[
R = \left[\begin{matrix}
1 & \ldots & 1 \\
\lambda_1(\pi_1) - \lambda_1(\pi_2) & \ldots & 	\lambda_r(\pi_1) - \lambda_r(\pi_2) 
\end{matrix}\right].
\] 
We show that $\rho(R) = 2$, where $\rho(R)$ denotes the rank of matrix $R$.  If this is the case, then $\dim(\Delta_r(w,\ell))= r - \rho(R) = r-2$. By contradiction, if $\rho(R) = 1$, there are two possibilities:
 \begin{itemize}
 	\item the second row is a multiple of the first: in this case the system has no solution;
 	\item the second row consists only of zeros, that is, $\lambda_h(\pi)=\lambda_h(\pi')$ for all $h\in \{1,\ldots,r\}$: this is not allowed by Assumption~\ref{ass:dist}.
 \end{itemize}
 Since $\dim(\Delta_r(w,\ell))< \dim(\Delta_r)$, the subspace of weights such that there exist at least two paths with same cost has null measure in $\Delta_r$. 
 \end{proof}
 In view of this lemma, we can prove the following proposition that basically states that function $\bar{g}$ is piecewise constant.
\begin{proposition}
\label{prop:piece}
For all $\bar{p} = ( \bar{p}^1,\ldots,\bar{p}^q)\in \Delta_r^q$ except those over a set of null measure over $\Delta_r^q$, it holds that there exists $\delta >0$ such that 
\[
(\forall p \in I_{\delta}(\bar{p}))\ \bar{g} (p) = \bar{g} (\bar{p}),
\]
where
\[
I_{\delta}(\bar{p}) = \{ p  = ( p^1,\ldots,p^q)\in \Delta_r^q\  | \;  \| p^\ell-\bar{p}^\ell\|< \delta,\ \ell \in \{1, \ldots, q\}\},
\]
 is a neighborhood of $\bar{p}$.
\end{proposition}   
\begin{proof}
Let $\bar{W}=W\setminus\left[\cup_{w\in W}\cup_{\ell=1}^q \Delta_r(w,\ell)\right]$.
For each O-D pair $w \in W$ and each $\ell \in \{1,\ldots,q\}$, let $\bar{p}^\ell \not \in \Delta_r(w,\ell)$, that is, $(\bar{p}^1,\ldots,\bar{p}^q)\in \bar{W}$. Note that, in view of Lemma~\ref{lem:null}, $\cup_{w\in W}\cup_{\ell=1}^q \Delta_r(w,\ell)$ is a set of null measure in $\Delta_r^q$. For each $w\in \bar{W}$, there exists a unique SP $\pi^{w,\ell}\in S^{w,\ell}$ for all $\ell \in \{1,\ldots,q\}$.
This means that there exist $a_{ij}^{w,\ell}, b_{ij}^{w,\ell}$ with:
\[
(\forall (i,j) \in A)\ (\forall \ell \in \{1,\ldots,q\})\ a_{ij}^{w,\ell}<c_{ij}^\top \bar{p}^\ell < b_{ij}^{w,\ell},
\]
such that if all arcs $(i,j)\in A$ have costs lying in $[a_{ij}^{w,\ell}, b_{ij}^{w,\ell}]$, then the SP $\pi^{w,\ell}$ remains the same as the one with costs of the arcs equal to $c_{ij}^\top \bar{p}^\ell$.
Now, let 
\[ [a_{ij}^\ell, b_{ij}^\ell ]= \bigcap_{w \in W} [a_{ij}^{w,\ell}, b_{ij}^{w,\ell}].\]
Note that
$a_{ij}^\ell<c_{ij}^\top \bar{p}^\ell < b_{ij}^\ell$ for all $(i,j)\in A$ and $\ell\in\{1,\ldots,q\}$,
must also hold.
Then, if for some $(p^1,\ldots,p^q)$ it holds that:
\begin{equation}		
\label{inequalities}
(\forall (i,j)\in A)\ (\forall \ell\in\{1,\ldots,q\})\ a_{ij}^\ell \leq c_{ij}^\top p^\ell \leq b_{ij}^\ell,	
\end{equation}
then $\bar{g}(p^1,\ldots,p^q)=\bar{g}(\bar{p}^1,\ldots,\bar{p}^q)$.
Now, let $\beta_{ij}^\ell = \min\{ c_{ij}^\top \bar{p}^\ell- a_{ij}^\ell, b_{ij}^\ell - c_{ij}^\top \bar{p}^\ell\}>0 $  for each arc $(i,j)\in A$, and $\beta^\ell= \min_{(i,j)\in A} \{\beta_{ij}^\ell\}>0$. 
Moreover, after setting $c= \max_{(i,j)\in A} \{ \|c_{ij}\|\}>0$, for all $\ell \in \{1,\ldots,q\}$ we define $\delta_\ell= \frac{\beta^\ell}{c}>0$. It follows that if $\| p^\ell-\bar{p}^\ell\|< \delta$, then for all $(i,j)\in A$:
\[| c_{ij}^\top (p^\ell-\bar{p}^\ell) | \leq \| p^\ell-\bar{p}^\ell\|  c  \leq \beta^\ell,\]
from which the inequalities~\eqref{inequalities} hold. Therefore, the thesis is proved with $\delta := \min_{\ell\in \{1,\ldots,q\}} \{\delta_\ell\}$.
\end{proof}
As a consequence of Proposition~\ref{prop:piece}, minimization of function $\bar{g}$ cannot be performed through a gradient search, since $\bar{g}$ is not even continuous. A combinatorial search through the constant pieces of $\bar{g}$ has to be performed.

\section{Data collection}
\label{sec:data}
Before discussing a possible approach to tackle the problem of identifying an unknown set of weights and the related probabilities, we point out that the definition of the problems discussed in Sections~\ref{sec:known} and~\ref{sec:unknown} requires the knowledge of some data. Namely, we need:
\begin{itemize}
\item vectors $c_1,\ldots,c_r$ with the basic costs;
\item set $W$ of O-D pairs with the related demand values (i.e., the number of users $u_w$ moving between origin $o_w$ and destination $d_w$, for each $w\in W$);
\item flows $\bar{x}_{ij}$ along (a subset of) the arcs of the network (i.e., the observed number of users traveling along the arcs).
\end{itemize}
In our experiments we will use synthetic data since we are still collecting real data, in particular for an application on the bike network of the City of Parma, Italy.
But it is important to discuss where and how these data are collected. 

Concerning the basic costs $c_1,\ldots,c_r$, different features may be considered.
These include traveling distance (or traveling time), safety, environmental conditions, such as pollution, temperature, pollens, or noise.
Note that some of these costs can be easily computed.
For instance, traveling distance can be obtained, e.g., from {\tt OpenStreetMap}.
Instead, other features, such as safety, are more qualitative and some way to convert them into quantities is needed.
To this end, also a collaboration with researchers more interested in the structural and infrastructural aspects of bike lanes might be fruitful.
It is also important to notice that some features have a static value (again, traveling distance), but some others have a dynamic value, depending on the time of the day or on the season of the year.
For instance, safety of a poorly enlightened bike lane decreases at night, while temperature or pollens are, of course, seasonal features.  

Concerning O-D pairs with related demand values, and observed flows, we have different ways to evaluate them. 
The measurement of flows of bikes within a town is possible through cameras or inductive loops, which are able to register bike passages. 
The advantage of such instruments is that they collect data referred to the whole population of bike users.
On the other hand, the limitation is that the population of bike users and, most of all, their O-D pairs with related demands (so called O-D matrix) might be rather difficult to estimate.
Therefore, alternative ways to collect data need to be explored.
Under this respect, a good source of information is represented by bike-sharing data.
Bike-sharing services are able to track the paths followed by users from a given starting point to a given ending point, thus providing information both about O-D pairs and about the travelled arcs, which is exactly what we need.
Note that we can include in the collection not only data referred to bikes but also to e-bikes and scooters, since habits of users of these means of transport are rather similar.
Such collection of data refers to a more limited population with respect to data coming from cameras, since only bike-sharing users are monitored, thus excluding all users with their own bike (and possibly the older segment of the general population).
Thus, inference from this data is less reliable with respect to data coming from cameras.
However, a possible way to limit this drawback is to include in the collection other data coming from different sources.
For instance, some public and private companies encourage the use of means of transport alternative to cars to reach the workplace.
To this end, they ask the employees to install apps through which it is possible to track their paths from their home to the workplace.
These are precious data for the definition of our problems and are a possible way to enlarge the sample of users.

We make a final remark about the collected data.
No matter how these data are collected, some care is needed when using them.
Indeed, here we are assuming that users are excellent optimizers: they have a cost function in mind, obtained as a convex combination of the basic costs, and they follow the best path with respect to that cost function.
In practice, this is not always the case.
Some users may make mistakes and follow suboptimal paths.
However, data of such users can be included in our model since they act like a noise signal in the measured traffic flow. Instead, some other users are simply roaming around without having any particular objective in mind.
These users should be considered as outliers and the data related to them should be removed from the definition of the problem.
In this paper we will not address the problem of detecting outliers, but this is certainly one of the problems to be addressed in future works.   

Finally, one limitation of this approach is that it assumes that all users make their choice based solely on the objective functions that we considered.
In order to make the model more and more realistic, we should try and take into account a growing number of factors that may influence route choices in cycling networks.
\section{Algorithm}
\label{sec:algo}
The idea of the proposed algorithm is to exploit the properties 
of objective function $g$ proved in Section~\ref{sec:unknown}  to generate a sequence of collections of weights with decreasing values of $g$.
The collection of weights is initialized with a sparse grid of weights $P_1$ over $\Delta_r$ (line~\ref{firstPartBegin}), for which $g(P_1)$ and the associated vector of probabilities $\alpha^{P_1}$ are computed by running procedure $\textsc{Identification}(P_1)$ (line~\ref{compgP1}), as described in Section~\ref{sec:known}.
Next, we enter the loop at lines~\ref{initrepeat}--\ref{firstPartEnd}. 
At each iteration of the loop, taking into account Corollary~\ref{coroll:1}, at line~\ref{c1} we restrict the attention to collection $\bar{P}$ of weights in $P_t$ whose probability is larger than a threshold $\epsilon$.
Next, taking into account Proposition~\ref{prop:monot}, at line~\ref{c2} we add new weights $P^{new}$ by perturbing those in $\bar{P}$ through procedure $\textsc{Perturbation}(\bar{P},\frac{1}{2^t})$, where $\frac{1}{2^t}$ is the size of the perturbation (thus decreasing with the iteration counter $t$). 
We do not detail procedure $\textsc{Perturbation}$, since different implementations are possible. 
For instance, if $r=3$ (the case that will be discussed in Section~\ref{sec:compexp}) we set:
\[
P^{new} = \left\{p+(i,j,-i-j)\frac{1}{2^t}\in \Delta_3 \mid
\substack{\displaystyle p \in \bar{P},\  i,j \in \{-1,0,1\},\\ \displaystyle i \neq 0 \vee j \neq 0}\right\}.
\]
At line~\ref{c3} we discard from further consideration the subset of weights $P^{del}\subset P_t \setminus \bar{P}$, whose probability is negligible and  whose distance from the closest weight in $\bar{P}$ is
at least $\frac{1}{2^{t -1}}$. The new collection $P_{t+1}$ is then defined at line~\ref{c4}, while at line~\ref{c5} we compute the new vector of probabilities $\alpha^{P_{t+1}}$ and value $g(P_{t+1})$. 
At line~\ref{c6} we increment by one the iteration counter and we double the value of threshold $\epsilon$. The loop is repeated until $\epsilon$ is lower than a predefined tolerance value $tol_1$ and the value of $g$ 
is decreased by at least a fraction $tol_2$ with respect to the previous iteration (see line~\ref{firstPartEnd}).
Once we exit the loop, at line~\ref{secondPartBegin} we call procedure $\textsc{Find Clusters}(P_t)$, which returns the set $B$ of the barycenters of the clusters identified within the collection of weights $P_t$, together with
vector $radius$, whose entries are the radii of the identified clusters. 
Note that the barycenters are computed by taking into account both the weights and their probability, that is, for a given cluster $P_c \subset P_t$, its barycenter is $b_c=\sum_{p\in P_c} \alpha^{P_t}_p p$.
The collection of weights is thus reduced from $P_t$ to $B$, and $B$ is employed to initialize the reduced collection
of  weights $P^*$ (line~\ref{b1}). In the final part of the algorithm, taking into account Proposition~\ref{prop:piece} and the need for a combinatorial search in the neighborhood of the weights, we try to refine each member of $P^*$ through a local search
({\tt For} loop at lines~\ref{b2}--\ref{secondPartEnd}).
For each $p^\ell\in P^*$ a search radius $\rho$ is initialized with the radius of the corresponding cluster (line~\ref{b3}). Next, the {\tt While} loop at lines~\ref{b4}--\ref{b8} is repeated. At line~\ref{b5} of the loop a local search around $p^\ell$ is performed with size of the perturbation equal to $\rho$. The local search explores the neighborhood $N(p^\ell)$ of $p^\ell$. If a new weight $p\in N(p^\ell)$ is identified such that 
$g( P^*\cup \{p\} \setminus \{p^\ell\})< g(P^*)$, then $P^*$ is updated into $P^*\cup \{p\} \setminus \{p^\ell\}$. Otherwise, the search radius $\rho$ is halved.

The loop is stopped as soon as the search radius falls below threshold $tol_3$.
Indeed, as seen in Proposition~\ref{prop:piece}, in a sufficiently small neighborhood of $P^*$, $g$ is constant. Thus, it makes sense to impose a lower limit for the search radius.
\alglanguage{pseudocode}
\begin{algorithm}
	\begin{algorithmic}[1]
		\Require $G = (V, A)$, $W$, $[\bar{x}_{ij}, (i,j) \in \bar{A}]$, $[u^w, w \in W]$, $tol_1,tol_2,tol_3\in [0,1]$.
		\Ensure $\alpha^{P^*}$, $P^{*}$.
		
		\State Initialize $P_1$, $\epsilon \gets 10^{-5}$, $t\gets 1$\label{firstPartBegin}
		\State{$[\alpha^{P_1}, g(P_1))] \gets \textsc{Identification}(P_1)$} \label{compgP1}
		\Repeat\label{initrepeat} 
		\State $\bar{P} \gets \{p \in P_t \mid \alpha^{P_t}_p > \epsilon\}$ \label{c1}
		\State $P^{new} \gets \textsc{Perturbation}(\bar{P},\frac{1}{2^{t}})$ \label{c2}
		\State $P^{del}$ is the set of weights in $P_t \setminus \bar{P}$ such that the minimum distance from $\bar{P}$ is greater than $\frac{1}{2^{t-1}}$. \label{c3}
		\State $P_{t+1} \gets P_{t} \cup P^{new} \setminus P^{del}$ \label{c4}
		\State{$[\alpha^{P_{t+1}}, g(P_{t+1}))] \gets \textsc{Identification}(P_{t+1})$} \label{c5}
		\State $\epsilon \gets 2\epsilon $, $t \gets t +1$ \label{c6}
	          \Until{$\epsilon<tol_1 \text{ and }  g(P_{t+1})\leq tol_2 \cdot g(P_t)$}\label{firstPartEnd}
		\State $[B, radius] \gets \textsc{Find Clusters}(P_t)$\label{secondPartBegin}
		\State $P^* \gets B$ \label{b1}
		\For{$p^\ell \in P^*$} \label{b2}
		\State $\rho \gets radius(\ell)$ \label{b3}
		\While {$\rho > tol_3$} \label{b4}
		\State $[P_{out},\rho_{out}] \gets \textsc{Local Search}(P^*, p^\ell, \rho)$ \label{b5}
		\State	$P^*\gets P_{out}$, $\rho\gets \rho_{out}$ 
		\EndWhile \label{b8}
		\EndFor\label{secondPartEnd}
		
	\end{algorithmic}
	\caption{Cyclists' route choice identification algorithm}
	\label{alg:choiceIdentification}
\end{algorithm}

\begin{algorithm}
	\begin{algorithmic}[1]
		\Require $P^*$, $p^\ell$, $\rho$
		\Ensure $P_{out}$, $\rho_{out}$
		\State $N(p^\ell) \gets \textsc{Perturbation}(\{p^\ell\}, \rho)$
		\State $P_{out}\gets P^*$, $\rho_{out}\gets \rho$
		\For {$p \in N(p^\ell)$}
		\State	$\bar{P} \gets P^* \cup \{p \}\setminus \{p^\ell\}$
		\If{$g(\bar{P})<g(P_{out})$}
		\State $P_{out} \gets \bar{P}$
		\EndIf
		\EndFor
		\If {$P_{out} = P^*$} 
		\State $\rho_{out}\gets \rho/2$
		\EndIf
	\end{algorithmic}
	\caption{\textsc{Local Search}}
	\label{alg:localsearch}
\end{algorithm}
Now we illustrate how the algorithm works through an example.
\begin{xmpl}
	We consider a set of five real weights $P_{ref}\in\Delta_3$ and we initialize $P_1$ with a grid of six equally distributed weights. In the top-left of Figure~\ref{fig:AA}, we show
	the projection of $\Delta_3$ onto $\mathbb{R}^2$. The stars are the real weights, while the empty and full circles are the weights in $P_1$. Function \textsc{Identification}($P_{1}$)  returns  probabilities $\alpha^{P_1}$. 
	Empty circles are weights in $P_1$ with a probability lower than $\epsilon$, while all other weights in $P_1$ are represented by full circles.
At some iteration, the resulting set $P_{t}$ is shown in the top-right of Figure~\ref{fig:AA}. Note that the weight distribution is denser with respect to $P_1$, and weights with higher probability are close to the real weights. Moreover, separate areas, which in the end will result in distinct clusters, are starting to appear. The search is concentrating in the neighborhood of real weights.
Then, when we exit the first {\tt While} loop, the weights in $P_t$ are divided into clusters based on their location. For each of them, the barycentre is found. In the bottom-left of Figure~\ref{fig:AA} we show one of the clusters with the corresponding barycentre (the square). The barycentre will replace the whole cluster of weights. Note that since a single weight replaces a set of weights, in this phase function $g$ may increase. But what we observed is that the increase is rather mild. 
Finally, for each barycentre, we perform a local search to refine the solution. In the bottom-right of Figure~\ref{fig:AA} we show the barycentre of the previous cluster and the corresponding final weight (the diamond) returned by the local search.

\begin{figure}[!htb]
\centering
\includegraphics[width=0.5\textwidth, clip=true, trim=27 21 33 18]{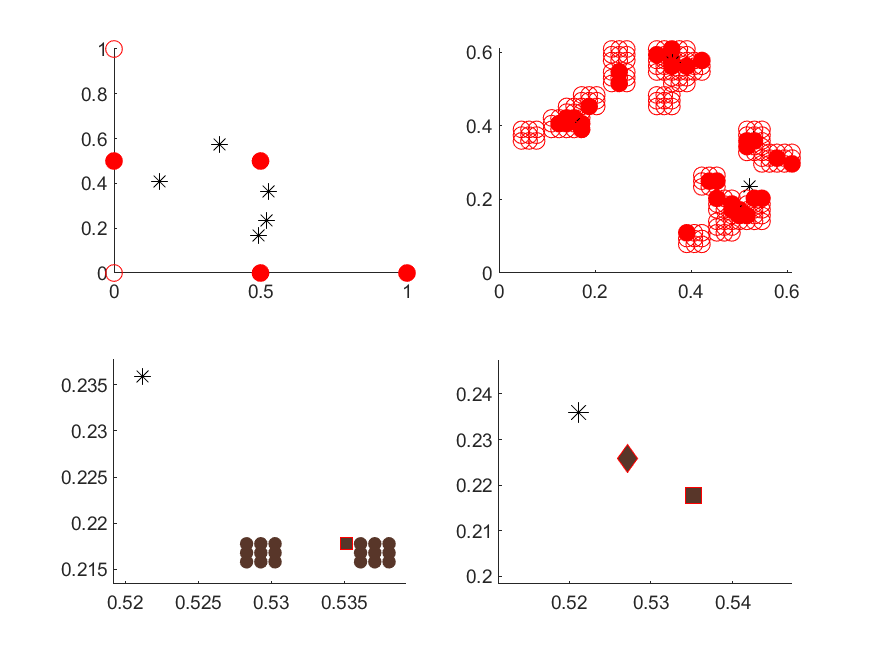}
\caption{Evolution of Algorithm~\ref{alg:choiceIdentification} over an instance.}
\label{fig:AA}           
\end{figure}
\end{xmpl}

\section{Simulations}\label{sec:compexp}
\begin{figure}[!htb]                                                                                                                           ~
	\begin{center} 
		\includegraphics[width=0.5\textwidth, clip=true, trim = 50 0 30 0]{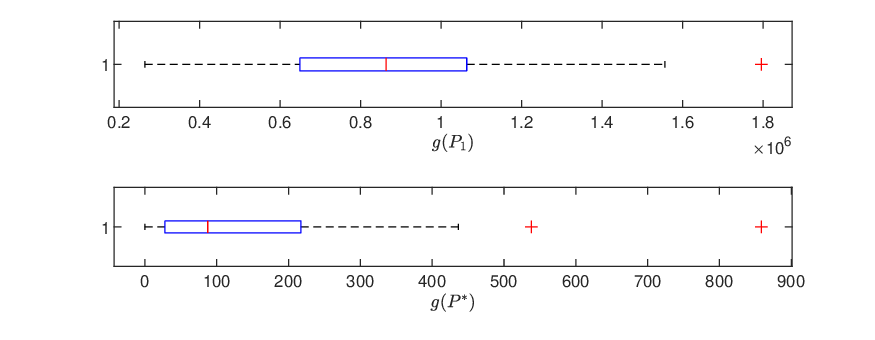}

		\caption{Distribution of the initial and final values of the objective function $g$.}
\label{fig:6}

\end{center}
\end{figure}

In order to test Algorithm~\ref{alg:choiceIdentification}, we perform some experiments on synthetic data. 
We consider the case of three distinct basic cost functions $c^1,c^2,c^3$ (i.e., $r=3$). These costs may represent, e.g., distance, safety and environmental conditions costs. 
The costs are integers randomly generated in $[5,20]$. They are normalized in such a way that   
$\sum_{(i,j)\in A}c^1_{ij}= \sum_{(i,j)\in A}c^2_{ij}=\sum_{(i,j)\in A}c^3_{ij}$.
We generate 100 instances for which the optimal set of weights $P_{ref}$ and the related probabilities $\alpha^{P_{ref}}$ are known in advance.
We randomly generate five distinct weights in $P_{ref}$ in such a way that the distance between them is at least $0.05$.
We also randomly generate the probabilities $\alpha^{P_{ref}}$ in such a way that all of them are not lower than 0.05.
We consider a grid graph $G=(V,A)$ with 1600 nodes, 
and a set $W$ with 1,000 O-D pairs. For each $w\in W$, we set $u_w=10$.
Next, for each $(i,j)\in A$ we calculate $\bar{x}_{ij}$ through~\eqref{eq:fij}--\eqref{eq:compxij} with $P=P_{ref}$ and $\alpha=\alpha_{ref}$, 
and we randomly select a subset of arcs $\bar{A}\subset A$, such that $|\bar{A}|=0.4|A|$. Note that this way
$[\alpha_{ref}, g(P_{ref})]=\textsc{Identification}(P_{ref})$, with $g(P_{ref})=0$. Finally, we run Algorithm~\ref{alg:choiceIdentification} after setting $tol_1=0.01$, $tol_2=0.85$, and $tol_3=0.005$.
All experiments have been performed on an Intel\textsuperscript{\textregistered} Core\textsuperscript{\texttrademark} i7-4510U CPU @ 2.60 GHz processor with 16 GB of RAM. Algorithm~\ref{alg:choiceIdentification} has been implemented in {\tt Matlab}. Shortest path problems have been solved by the {\tt Matlab} routine implementing Dijkstra's algorithm.
The convex QPs have been solved through {\tt Gurobi} called via {\tt Yalmip}. Clusters have been identified through the {\tt Matlab} routine for Hierarchical Clustering.
In Figure~\ref{fig:6} we show the distribution of the initial value $g(P_1)$ (upper picture) and of the final value $g(P^*)$ (lower picture) for all the tested instances. Note that the final value is significatively lower than the initial value $g(P_1)$, which is larger than $10^5$.
The large part of the improvement with respect to $g(P_1)$ is due to the first {\tt While} loop. As previously commented, when we move from the set $P_t$ returned by this loop to the set $B$ of barycenters, there is a small increase of $g$
(but with the advantage of having significatively reduced the number of weights). The final local search is able to refine the set of barycenters and allows for a further mild reduction of $g$.

\begin{figure}[!htb]
	\begin{center} 
		\includegraphics[width=0.5\textwidth, clip=true, trim = 40 0 30 0]{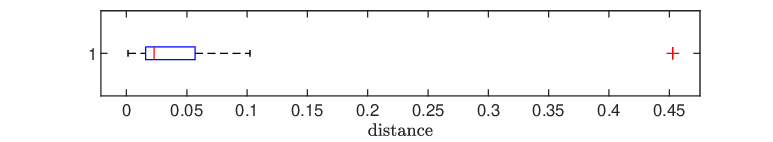}
		\caption{Distribution of the distances between the real weights and the weights returned by Algorithm~\ref{alg:choiceIdentification}.}
		\label{fig:9}
	\end{center}
\end{figure}

\begin{figure}[!htb]
	\begin{center} 
		\includegraphics[width=0.5\textwidth, clip=true, trim = 50 0 30 0]{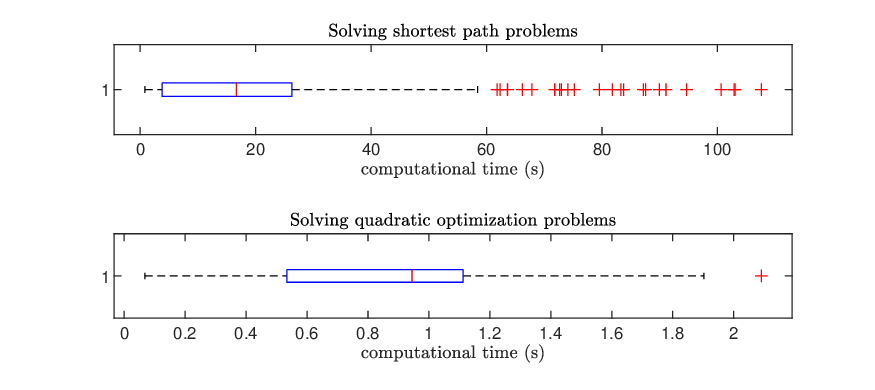}
		\caption{Computational time to solve the SP problems compared to the computational time to solve the convex QPs.}
		\label{fig:10}
	\end{center}
\end{figure}

In Figure~\ref{fig:9} we show the distribution of the Euclidean distances between the real weights and the weights $P^*$ returned by Algorithm~\ref{alg:choiceIdentification}. As we can see, the distances are, with a single exception, rather small.
The main computational cost in Algorithm~\ref{alg:choiceIdentification} is represented by the calls of the procedure $\textsc{Identification}$. In turn, the computing times of these calls are determined by: (i) the solution
of SP problems; (ii) the solution of convex QPs. Then, in Figure~\ref{fig:10} we compare the cumulative time needed
by the solutions of the former (upper picture) and by the latter (lower picture). It appears quite clearly that computing times are dominated by the solutions of the SP problems.

\section{Conclusions and future work}
In this paper we proposed an optimization model and a related solution algorithm for the identification of the criteria with which different users select their paths when moving
in a bike network. We assume that users move along SPs but that the costs of these paths are convex combinations of a few basic costs, which take into account different aspects, such as distance or safety. 
The proposed optimization model is based on the observation of real flows of users within the network. First, a simplified problem with a set of weights (coefficients of the convex combinations giving the costs optimized by the users) known in advance is tackled. Such problem is solved through a polynomial-time algorithm based on the solution of many SP problems and a single convex QP problem. Next, an algorithm to identify an unknown set of weights is proposed.
Experiments over synthetic data are reported. In a future work, the proposed methodology will be applied to the real case of the bike network of Parma, Italy, using the data made available by bike sharing services.
This will require a careful identification and definition of the basic costs. Moreover, it will be necessary to identify outliers (i.e., users who are simply biking around without any optimization cost in mind).

\bibliographystyle{abbrv}
\bibliography{biblio}

\end{document}